\documentclass[10pt]{article}
\usepackage{mathrsfs}
\usepackage{amsthm}
\usepackage{amssymb}
\usepackage{amsmath}
\usepackage{graphicx}
\usepackage{color}
\usepackage{amsfonts}
\usepackage{float}
\usepackage{cite}
\usepackage[text={140mm,210mm},left=45mm,vmarginratio=1:1]{geometry}
\newtheorem{theorem}{Theorem}[section]

\newtheorem{lemma}[theorem]{Lemma}

\newtheorem{corollary}[theorem]{Corollary}

\numberwithin{equation}{section}
\normalsize

\begin{document}
\title{\textbf{An improved upper bound for critical value of the contact process on $\mathbb{Z}^d$ with $d\geq 3$}}

\author{Xiaofeng Xue \thanks{\textbf{E-mail}: xfxue@bjtu.edu.cn \textbf{Address}: School of Science, Beijing Jiaotong University, Beijing 100044, China.}\\ Beijing Jiaotong University}

\date{}
\maketitle

\noindent {\bf Abstract:} In this paper we give an improved upper
bound for critical value $\lambda_c$ of the basic contact process on
the lattice $\mathbb{Z}^d$ with $d\geq 3$. As a direct corollary of
out result,
\[
\lambda_c\leq 0.340
\]
when $d=3$.

\quad


\section{Introduction}\label{section one}

In this paper we are concerned with the basic contact process on
$\mathbb{Z}^d$ with $d\geq 3$. First we introduce some notations.
For each $x=(x_1,x_2,\ldots,x_d)\in \mathbb{Z}^d$, we use $\|x\|$ to
denote the $l_1$-norm of $x$, i.e.,
\[
\|x\|=\sum_{i=1}^d|x_i|.
\]
For any $x,y\in \mathbb{Z}^d$, we write $x\sim y$ when end only when
$\|x-y\|=1$, i.e., $x\sim y$ means that $x$ and $y$ are neighbors on
$\mathbb{Z}^d$. For $1\leq i\leq d$, we use $e_i$ to denote the
$i$th elementary unit vector of $\mathbb{Z}^d$, i.e.,
\begin{equation}\label{equ 1.1}
e_i=(0,\ldots,0,\mathop 1\limits_{i \text{th}},0,\ldots,0).
\end{equation}
We use $O$ to denote the origin of $\mathbb{Z}^d$.

The contact process $\{\eta_t\}_{t \geq 0}$ on $\mathbb{Z}^d$ is a
spin system with state space $\{0,1\}^{\mathbb{Z}^d}$ (see the
definition of the spin system in Chapter 3 of \cite{Lig1985}). The
flip rates function of $\{\eta_t\}_{t\geq 0}$ is given by
\begin{equation}\label{equ 1.2 transition rate}
c(x,\eta)=
\begin{cases}
1 & \text{~if~} \eta(x)=1,\\
\lambda\sum_{y:y\sim x}\eta(y) & \text{~if~} \eta(x)=0
\end{cases}
\end{equation}
for any $(\eta,x)\in \{0,1\}^{\mathbb{Z}^d}\times \mathbb{Z}^d$,
where $\lambda>0$ is a constant called the infection rate. That is
to say, the state of the process flips from $\eta$ to $\eta^x$ at
rate $c(x,\eta)$, where
\[
\eta^x(y)=
\begin{cases}
\eta(y) & \text{~if~}y\neq x,\\
1-\eta(x) & \text{~if~}y=x.
\end{cases}
\]

Intuitively, the contact process describes the spread of an epidemic
on the graph. Vertices in state $1$ are infected while that in state
$0$ are healthy. An infected vertex waits for an exponential time
with rate $1$ to become healthy while an healthy one is infected at
rate proportional to the number of infected neighbors.

The contact process is introduced by Harris in \cite{Har1974}. For a
detailed survey of the study of the contact process, see Chapter 6
of \cite{Lig1985} and Part one of \cite{Lig1999}.

In this paper we are mainly concerned with the critical value of the
contact process. Assuming that $\eta_0(x)=1$ for any $x\in
\mathbb{Z}^d$, then the critical value $\lambda_c$ is defined as
\begin{equation}\label{equ 1.3 definition of critical value}
\lambda_c=\sup\big\{\lambda:\lim_{t\rightarrow+\infty}P_\lambda(\eta_t(O)=1)=0\big\},
\end{equation}
where $P_\lambda$ is the probability measure of the contact process
with infection rate $\lambda$. The definition of $\lambda_c$ is
reasonable according to the following property of the contact
process. For $\lambda_1\geq \lambda_2$ and $t>s$, conditioned on all
the vertices are in state $1$ at $t=0$,
\begin{equation}\label{equ 1.4}
P_{\lambda_1}(\eta_s(O)=1)\geq P_{\lambda_2}(\eta_t(O)=1).
\end{equation}
A rigorous proof of Equation \eqref{equ 1.4} is given in Section 6.1
of \cite{Lig1985}.

When $d=1$, it is shown in Section 6.1 of \cite{Lig1985} that
$\lambda_c(1)\leq 2$. Liggett improves this result in \cite{Lig1995}
by showing that $\lambda_c(1)\leq 1.94$. For $d\geq 3$, it is shown
in \cite{Hol1981} that
\[
\lambda_c(d)\leq \alpha_1(d)=\frac{1}{\gamma_d}-1
\]
while it is shown in \cite{Gri1983} that
\[
\lambda_c(d) \leq \alpha_2(d)=\frac{1}{2d(2\gamma_d-1)},
\]
where $\gamma(d)>1/2$ is the probability that the simple random walk
on $\mathbb{Z}^d$ starting at $O$ never returns to $O$. Both these
two results lead to the conclusion that
\[
\lim_{d\rightarrow+\infty}2d\lambda_c(d)=1.
\]
When $d=3$, according to the well-known result that $\gamma_3\approx
0.659$,
\[
\alpha_1(3)=0.517<\alpha_2(3)=0.523.
\]
However, $\alpha_2(d)<\alpha_3(d)$ for sufficiently large $d$
according to the fact that
\[
\frac{1}{\gamma_d}-1=\frac{1}{2d}+\frac{3}{4d^2}+o(\frac{1}{d^2})
\]
while
\[
\frac{1}{2d(2\gamma_d-1)}=\frac{1}{2d}+\frac{1}{2d^2}+o(\frac{1}{d^2}).
\]
In this paper, we will give another upper bound $\beta(d)$ for the
critical value $\lambda_c(d)$ when $d\geq 3$. $\beta(d)$ satisfies
that $\beta(d)<\min\{\alpha_1(d), \alpha_2(d)\}$ for each $d\geq 3$. For the precise result, see the next
section.

\section{Main result}\label{section two}
In this section we will give our main result. First we introduce
some notations and definitions. From now on we assume that at $t=0$
all the vertices on $\mathbb{Z}^d$ are in state $1$ for the contact
process, then let $\lambda_c$ be the critical value of the contact
process defined as in Equation \eqref{equ 1.3 definition of critical
value}. We write $\lambda_c$ as $\lambda_c(d)$ when we need to point
out the dimension $d$ of the lattice. We denote by $\{S_n\}_{n\geq
0}$ the simple random walk on $\mathbb{Z}^d$, i.e.,
\[
P\big(S_{n+1}=y\big|S_n=x\big)=\frac{1}{2d}
\]
for each $y$ that $y\sim x$ and $n\geq 0$. We define
\[
\gamma=P\big(S_n\neq O\text{~for all~}n\geq 1\big|S_0=O\big)
\]
as the probability that the simple random walk never return to $O$
conditioned on $S_0=O$. We write $\gamma$ as $\gamma_d$ when we need
to point out the dimension $d$ of the lattice.

The following theorem gives an upper bound of $\lambda_c(d)$ for
$d\geq 3$, which is our main result.

\begin{theorem}\label{theorem 2.1 main}
For each $d\geq 3$,
\[
\lambda_c(d)\leq \frac{2-\gamma_d}{2d\gamma_d}.
\]
\end{theorem}

It is shown in \cite{Gri1983} that $\lambda_c(d)\leq
\alpha_2(d)=\frac{1}{2d(2\gamma_d-1)}$ for each $d\geq 3$. Since
$\gamma_d<1$,
\[
(2-\gamma_d)(2\gamma_d-1)-\gamma_d=-2(\gamma_d-1)^2<0
\]
and hence $\frac{2-\gamma_d}{2d\gamma_d}<\alpha_2(d)$ for each $d\geq 3$.
It is shown in \cite{Hol1981} that $\lambda_c(d)\leq
\alpha_1(d)=\frac{1}{\gamma_d}-1$ for each $d\geq 3$. By direct calculation,
\begin{align*}
1-\gamma &\geq P\big(S_2=O\big|S_0=O\big)+P\big(S_4=O,S_2\neq O\big|S_0=O\big)\\
&=\frac{4d^2+4d-3}{8d^3}>\frac{1}{2d-1}
\end{align*}
when $d\geq 3$ and hence $\frac{2-\gamma_d}{2d\gamma_d}<\alpha_1(d)$ for each $d\geq 3$.

For $d=3$, according to the well known result that $\gamma_3\approx
0.659$, we have the following direct corollary.
\begin{corollary}\label{corollary 2.2}
\[
\lambda_c(3)\leq \frac{2-\gamma_3}{6\gamma_3}\leq 0.340.
\]
\end{corollary}
This corollary improves the upper bound of $\lambda_c(3)$ given by
$\alpha_1(3)$, which is $0.517$. According to the example given in
Section 3.5 of \cite{Lig1985},
\[
\lambda_c(d)\geq \frac{1}{2d-1}
\]
for each $d\geq 1$ and hence $\lambda_c(3)\in [0.2,0.340]$.

We will prove Theorem \ref{theorem 2.1 main} in the next section. A
Markov process $\{\xi_t\}_{t\geq 0}$ with state space
$[0,+\infty)^{\mathbb{Z}^d}$ will be introduced as a main auxiliary
tool for the proof. The definition of $\{\xi_t\}_{t\geq 0}$ is
similar with that of the binary contact path process introduced in
\cite{Gri1983}, except for some modifications in several details.

\quad

\section{Proof of Theorem \ref{theorem 2.1 main}}\label{section
three} In this section we give the proof of Theorem \ref{theorem 2.1
main}. Throughout this section we assume that the dimension $d$ is
fixed and at least $3$, which ensures that $\gamma>\frac{1}{2}$. Our
aim is to prove the following lemma, Theorem \ref{theorem 2.1 main}
follows from which directly.
\begin{lemma}\label{Lemma 3.1}
If $a,b>0$ satisfies
\[
2(a+b-1)-(a^2+b^2-1)-2ab(1-\gamma)>0
\]
then
\[
\lambda_c\leq \frac{1}{2d\big(2(a+b-1)-(a^2+b^2-1)-2ab(1-\gamma)\big)}.
\]
\end{lemma}
If we choose $a=b=1$, then Lemma \ref{Lemma 3.1} gives the upper bound
of $\lambda_c$ the same as that given in \cite{Gri1983}. However,
the best choices of $a,b$ are $a=b=\frac{1}{2-\gamma}$, which gives the following proof
of Theorem \ref{theorem 2.1 main}.

\proof[Proof of Theorem \ref{theorem 2.1 main}] Let
$L(a,b)=2(a+b-1)-(a^2+b^2-1)-2ab(1-\gamma)$, then
\[
\sup\big\{L(a,b):a>0,b>0\big\}=L(\frac{1}{2-\gamma},\frac{1}{2-\gamma})=\frac{\gamma}{2-\gamma}.
\]
As a result, let $a=b=\frac{1}{2-\gamma}$, then
\[
\lambda_c\leq \frac{1}{2dL(a,b)}=\frac{2-\gamma}{2d\gamma}
\]
according to Lemma \ref{Lemma 3.1}.

\qed

The remainder of this paper is devoted to the proof of Lemma
\ref{Lemma 3.1}. From now on we assume that $a,b$ are positive
constants which satisfies
\[
2(a+b-1)-(a^2+b^2-1)-2ab(1-\gamma)>0.
\]
Let $\{\xi_t\}_{t\geq 0}$ be a continuous time Markov process with
state space $[0,+\infty)^{\mathbb{Z}^d}$ and generator function
given by
\begin{align}\label{equ 3.1}
\Omega f(\xi)=&\sum_{x\in
\mathbb{Z}^d}\big[f(\xi^{x,0})-f(\xi)\big]+\sum_{x\in
\mathbb{Z}^d}\sum_{y:y\sim
x}\lambda\big[f(\xi_{a,b}^{x,y})-f(\xi)\big]\\
&+\sum_{x\in
\mathbb{Z}^d}f_x^{\prime}(\xi)\Big(1-2d\lambda[(b-1)+a]\Big)\xi(x)\notag
\end{align}
for any $\xi\in [0,+\infty)^{\mathbb{Z}^d}$ and sufficiently smooth
function $f$ on $[0,+\infty)^{\mathbb{Z}^d}$, where
\begin{align*}
&\xi^{x,0}(y)=
\begin{cases}
\xi(y) & \text{~if~}y\neq x,\\
0 & \text{~if~}y=x,
\end{cases}\\
&\xi^{x,y}_{a,b}(z)=
\begin{cases}
\xi(z) & \text{~if~} z\neq x,\\
b\xi(x)+a\xi(y) & \text{~if~} z=x
\end{cases}
\end{align*}
and $f^{\prime}_x$ is the partial derivative of $f(\xi)$ with
respect to the coordinate $\xi(x)$.

If $a=b=1$, then $\{\xi_t\}_{t\geq 0}$ is the binary contact path
process introduced in \cite{Gri1983} after a time-scaling.
$\{\xi_t\}_{t\geq 0}$ belongs to a large crowd of continuous-time
Markov processes called linear systems. For the definition and basic
properties of the linear system, see Chapter 9 of \cite{Lig1985}.

According to the definition of $\Omega$, $\{\xi_t\}_{t\geq 0}$
evolves as follows. For each $x\in \mathbb{Z}^d$ and each neighbor
$y$ of $x$, $\xi_t(x)$ flips to $0$ at rate $1$ while flips to
$b\xi_t(x)+a\xi_t(y)$ at rate $\lambda$. Between the jumping moments
of $\{\xi_t(x)\}_{t\geq 0}$, $\xi_t(x)$ evolves according to the ODE
\begin{equation}\label{equ 3.2 ODE}
\frac{d}{dt}\xi_t(x)=\Big(1-2d\lambda\big[(b-1)+a\big]\Big)\xi_t(x).
\end{equation}
That is to say, if $\xi(x)$ does not jump during $[t,t+s]$, then
\[
\xi_{t+r}(x)=\xi_t(x)\exp\Big\{r\Big(1-2d\lambda\big[(b-1)+a\big]\Big)\Big\}
\]
for $0<r<s$.

The linear system $\{\xi_t\}_{t\geq 0}$ and the contact process
$\{\eta_t\}_{t\geq 0}$ have the following relationship.
\begin{lemma}\label{lemma 3.2}
For any $x\in \mathbb{Z}^d$ and $t\geq 0$, let
\[
\widehat{\eta}_t(x)=
\begin{cases}
1 & \text{~if~}\xi_t(x)>0,\\
0 & \text{~if~}\xi_t(x)=0,
\end{cases}
\]
then $\{\widehat{\eta}_t\}_{t\geq 0}$ is a version of the contact
process introduced in Equation \eqref{equ 1.2 transition rate}.
\end{lemma}

\proof[Proof of Lemma \ref{lemma 3.2}]

ODE \eqref{equ 3.2 ODE} can not make $\{\xi_t(x)\}_{t\geq 0}$ flip
from $0$ to a positive value or flip from a positive value to $0$,
hence $\widehat{\eta}_t(x)$ stays its value between jumping moments
of $\xi(x)$. If $\widehat{\eta}_t(x)=1$, i.e, $\xi_t(x)>0$, then
$\widehat{\eta}_t(x)$ flips to $0$ when and only when $\xi_t(x)$
flips to $0$ at some jumping moment. As a result,
$\widehat{\eta}_t(x)$ flips from $1$ to $0$ at rate $1$. If
$\widehat{\eta}_t(x)=0$, i.e, $\xi_t(x)=0$, then
$\widehat{\eta}_t(x)$ flips to $1$ when and only when $\xi_t(x)$
flips to
\[
b\xi_t(x)+a\xi_t(y)=a\xi_t(y)
\]
for a neighbor $y$ with $\xi_t(y)>0$ at some jumping moment. As a
result, $\widehat{\eta}_t(x)$ flips from $0$ to $1$ at rate
\[
\lambda\sum_{y:y\sim x}1_{\{\xi_t(y)>0\}}=\lambda\sum_{y:y\sim
x}\widehat{\eta}_t(y),
\]
where $1_A$ is the indicator function of the event $A$. In
conclusion, $\{\widehat{\eta}_t\}_{t\geq 0}$ evolves in the same way
as a contact process evolves according to the flip rates function
given in Equation \eqref{equ 1.2 transition rate}.

\qed

By Lemma \ref{lemma 3.2}, from now on we assume that
$\{\eta_t\}_{t\geq 0}$ and $\{\xi_t\}_{t\geq 0}$ are coupled under
the same probability space such that $\eta_0(x)=\xi_0(x)=1$ for each
$x\in \mathbb{Z}^d$ and $\eta_t(x)=1$ when and only when
$\xi_t(x)>0$.

The following two lemmas about expectations of $\xi_t(x)$ and
$\xi_t(x)\xi_t(y)$ are important for the proof of Lemma \ref{Lemma
3.1}.
\begin{lemma}\label{lemma 3.3}
If $\xi_0(x)=1$ for any $x\in \mathbb{Z}^d$, then
\[
E\xi_t(x)=1
\]
for any $x\in \mathbb{Z}^d$ and $t\geq 0$.
\end{lemma}

\begin{lemma}\label{lemma 3.4}
For any $x\in \mathbb{Z}^d$ and $t\geq 0$, let
$F_t(x)=E\big[\xi_t(O)\xi_t(x)\big]$, then conditioned on
$\xi_0(x)=1$ for all $x\in \mathbb{Z}^d$,
\begin{equation}\label{equ 3.2 two ODE}
\frac{d}{dt}F_t=\Big(\frac{d}{dt}F_t(x)\Big)_{x\in Z^d}=G_\lambda
F_t,
\end{equation}
where $G_\lambda$ is a $\mathbb{Z}^d\times \mathbb{Z}^d$ matrix that
\[
G_\lambda(x,y)=
\begin{cases}
-4a\lambda d & \text{~if~} x\neq 0 \text{~and~}x=y,\\
2a\lambda & \text{~if~} x\neq 0 \text{~and~}x\sim y,\\
1-4d\lambda(b-1)-4d\lambda a+2d\lambda(b^2-1)+2d\lambda a^2 & \text{~if~} x=y=0,\\
4abd\lambda &\text{~if~}x=0 \text{~and~}y=e_1,\\
0 &\text{~otherwise}
\end{cases}
\]
and $e_1$ is defined as in Equation \eqref{equ 1.1}.
\end{lemma}
Note that when we say $F_1=GF_2$ for functions $F_1,F_2$ on
$\mathbb{Z}^d$ and $\mathbb{Z}^d\times \mathbb{Z}^d$ matrix $G$, we
mean
\[
F_1(x)=\sum_{y\in Z^d}G(x,y)F_2(y)
\]
for each $x\in \mathbb{Z}^d$, as the product of finite-dimensional
matrixes.

The proofs of Lemmas \ref{lemma 3.3} and \ref{lemma 3.4} rely
heavily on Theorems 9.1.27 and 9.3.1 of \cite{Lig1985}. These two
theorems can be seen as the extension of the Hille-Yosida Theorem
for the linear system, which ensures that we can execute the
calculation
\begin{equation}\label{equ 3.3}
\frac{d}{dt}S(t)f=S(t)\Omega f
\end{equation}
for a linear system with generator $\Omega$ and semi-group
$\{S_t\}_{t\geq 0}$ when $f$ has the form $f(\xi)=\xi(x)$ or
$f(\xi)=\xi(x)\xi(y)$.

\proof[Proof of Lemma \ref{lemma 3.3}] By the generator $\Omega$ of
$\{\xi_t\}_{t\geq 0}$ and Theorem 9.1.27 of \cite{Lig1985} (i.e.,
Equation \eqref{equ 3.3} for $f(\xi)=\xi(x)$),
\[
\frac{d}{dt}E\xi_t(x)=-E\xi_t(x)+\lambda\sum_{y:y\sim
x}\Big[(b-1)E\xi_t(x)+aE\xi_t(y)\Big]+\Big(1-2d\lambda\big[(b-1)+a\big]\Big) E\xi_t(x)
\]
for each $x\in \mathbb{Z}^d$. Since $\xi_0(x)=1$ for all $x\in
\mathbb{Z}^d$, $E\xi_t(x)$ does not depend on the choice of $x$
according to the spatial homogeneity of $\{\xi_t\}_{t\geq 0}$.
Therefore,
\begin{align*}
\frac{d}{dt}E\xi_t(x)&=-E\xi_t(x)+\lambda\sum_{y:y\sim
x}\Big[(b-1)E\xi_t(x)+aE\xi_t(y)\Big]+\Big(1-2d\lambda[(b-1)+a]\Big) E\xi_t(x)\\
&=-E\xi_t(x)+2d\lambda (a+b-1)E\xi_t(x)+\big(1-2d\lambda(a+b-1)\big) E\xi_t(x)=0.
\end{align*}
As a result, $E\xi_t(x)\equiv E\xi_0(x)=1$.

\qed

\proof[Proof of Lemma \ref{lemma 3.4}] According to the generator
$\Omega$ of $\{\xi_t\}_{t\geq 0}$ and Theorem 9.3.1 of
\cite{Lig1985} (i.e., Equation \eqref{equ 3.3} for
$f(\xi)=\xi(x)\xi(y)$),
\begin{align}\label{equ 3.4}
\frac{d}{dt}F_t(x)=&-2F_t(x)+\lambda\sum_{y:y\sim
O}\Big((b-1)F_t(0)+aE\big[\xi_t(y)\xi_t(x)\big]\Big)\notag\\
&+\lambda\sum_{y:y\sim
x}\Big((b-1)F_t(0)+aF_t(y)\Big)
+2\Big(1-2d\lambda(a+b-1)\Big)F_t(x)
\end{align}
when $x\neq O$ while
\begin{align}\label{equ 3.5}
\frac{d}{dt}F_t(O)=&-F_t(O)+\lambda\sum_{y:y\sim
O}2abF_t(y)+2d\lambda(b^2-1)F_t(O)+\lambda\sum_{y:y\sim O}a^2E\big[\xi_t^2(y)\big]\notag\\
&+2\big(1-2d\lambda(a+b-1)\big)F_t(O).
\end{align}
Since $\xi_0(x)=1$ for any $x\in \mathbb{Z}^d$, according to the
spatial homogeneity of $\{\xi_t\}_{t\geq 0}$,
\[
E\big[\xi_t(x)\xi_t(y)\big]=F_t(y-x)=F_t(x-y)
\]
for any $x,y\in \mathbb{Z}^d$ and
\[
F_t(e_i)=F_t(-e_i)=F_t(e_1)
\]
for $1\leq i\leq d$. Therefore, by Equations \eqref{equ 3.4} and
\eqref{equ 3.5},
\begin{equation}\label{equ 3.6}
\frac{d}{dt}F_t(x)=
\begin{cases}
-4ad\lambda F_t(x)+2a\lambda\sum_{y:y\sim x}F_t(y) &
\text{~if~}x\neq O,\\
\big[1-4d\lambda(a+b-1)+2d\lambda(b^2-1)+2da^2\lambda\big]F_t(O)+4abd\lambda F_t(e_1) &
\text{~if~}x=O.
\end{cases}
\end{equation}
Lemma \ref{lemma 3.4} follows from Equation \eqref{equ 3.6}
directly.

\qed

The following lemma shows that if $\lambda$ ensures the existence of
an positive eigenvector of $G_\lambda$ with respect to the
eigenvalue $0$, then $\lambda$ is an upper bound of $\lambda_c$,
which is crucial for us to prove Lemma \ref{Lemma 3.1}.

\begin{lemma}\label{lemma 3.5}
If there exists $K: \mathbb{Z}^d\rightarrow [0,+\infty)$ that
$\inf_{x\in \mathbb{Z}^d}K(x)>0$ and
\[
G_\lambda K=0 ~(\text{here $0$ means the zero function on
$\mathbb{Z}^d$}),
\]
where $G_\lambda$ is defined as in Lemma \ref{lemma 3.4}, then
\[
\lambda\geq \lambda_c.
\]
\end{lemma}

We give the proof of Lemma \ref{lemma 3.5} at the end of this
section. Now we show how to utilize Lemma \ref{lemma 3.5} to prove
Lemma \ref{Lemma 3.1}.

\proof[Proof of Lemma \ref{Lemma 3.1}] Let $\{S_n\}_{n\geq 0}$ be
the simple random walk on $\mathbb{Z}^d$ as we have introduced in
Section \ref{section two}, then we define
\[
H(x)=P\big(S_n=O \text{~for some~}n\geq 0\big|S_0=x\big)
\]
for any $x\in \mathbb{Z}^d$. Then $H(O)=1$ and
\begin{equation}\label{equ 3.7}
H(x)=\frac{1}{2d}\sum_{y:y\sim x}H(y)
\end{equation}
for any $x\neq O$. According to the spatial homogeneity of the
simple random walk,
\begin{align}\label{equ 3.8}
\gamma&=P\big(S_n\neq O \text{~for all~}n\geq 1\big|S_0=O\big) \notag\\
&=P\big(S_n\neq O \text{~for all~}n\geq 0\big|S_0=e_1\big)=1-H(e_1).
\end{align}
For $a,b>0$ that
\[
2(a+b-1)-(a^2+b^2-1)-2ab(1-\gamma)>0
\]
and
$\lambda>\frac{1}{2d\big[2(a+b-1)-(a^2+b^2-1)-2ab(1-\gamma)\big]}$, we define
\[
K(x)=H(x)+\frac{2d\lambda\big[2(a+b-1)-(a^2+b^2-1)-2ab(1-\gamma)\big]-1}{1+2d\lambda(a+b-1)^2}
\]
for each $x\in \mathbb{Z}^d$. Then,
\[
\inf_{x\in \mathbb{Z}^d}K(x)\geq
\frac{2d\lambda\big[2(a+b-1)-(a^2+b^2-1)-2ab(1-\gamma)\big]-1}{1+2d\lambda(a+b-1)^2}>0
\]
and $G_\lambda K=0$ according to Equations \eqref{equ 3.7},
\eqref{equ 3.8} and the definition of $G_\lambda$. As a result, by
Lemma \ref{lemma 3.5},
\[
\lambda\geq \lambda_c
\]
for any $\lambda>\frac{1}{2d\big[2(a+b-1)-(a^2+b^2-1)-2ab(1-\gamma)\big]}$ and hence
\[
\lambda_c\leq \frac{1}{2d\big[2(a+b-1)-(a^2+b^2-1)-2ab(1-\gamma)\big]}.
\]

\qed

At last we give the proof of Lemma \ref{lemma 3.5}.

\proof[Proof of Lemma \ref{lemma 3.5}]

For any $x,y\in \mathbb{Z}^d$, we define
\[
G_\lambda^2(x,y)=\sum_{u\in
\mathbb{Z}^d}G_\lambda(x,u)G_\lambda(u,y).
\]
It is easy to check that the sum in the right-hand side converges
since only finite terms are not zero. By induction, if $G_\lambda^k$
is well-defined for $1\leq k\leq n$, then we define
\[
G_\lambda^{n+1}(x,y)=\sum_{u\in
\mathbb{Z}^d}G_\lambda^n(x,u)G_\lambda(u,y).
\]
It is easy to check that $G_\lambda^n$ is well-defined for each
$n\geq 1$ according to the definition of $G_\lambda$ and
\[
\sup_{x,y\in
\mathbb{Z}^d}\sum_{n=0}^{+\infty}\frac{t^n|G^n_\lambda(x,y)|}{n!}<+\infty
\]
for any $t\geq 0$, where $G_\lambda^0(x,y)=1_{\{x=y\}}$. Then, it is
reasonable to define the $\mathbb{Z}^d\times \mathbb{Z}^d$ matrix
$e^{tG_\lambda}$ as
\[
e^{tG_\lambda}(x,y)=\sum_{n=0}^{+\infty}\frac{t^nG_\lambda^n(x,y)}{n!}
\]
for $x,y\in \mathbb{Z}^d$ and $t\geq 0$. Since $K$ satisfies
$G_\lambda K=0$,
\[
G_\lambda^nK=G_\lambda^{n-1}G_\lambda K=0
\]
for each $n\geq 1$ and hence
\begin{equation}\label{equ 3.9}
(e^{tG_\lambda}K)(x)=\sum_{y\in
\mathbb{Z}^d}e^{tG_\lambda}(x,y)K(y)=\sum_{y\in
\mathbb{Z}^d}G^0_\lambda(x,y)K(y)=K(x)
\end{equation}
for each $x\in \mathbb{Z}^d$ and $t\geq 0$, i.e., $K$ is the
eigenvector of $e^{tG_\lambda}$ with respect to the eigenvalue $1$.

For any $\xi\in (-\infty,+\infty)^{\mathbb{Z}^d}$, we define
\[
\|\xi\|_{\infty}=\sup_{x\in \mathbb{Z}^d}|\xi(x)|.
\]
Furthermore, we define
\[
W=\{\xi\in (-\infty,+\infty)^{\mathbb{Z}^d}:
\|\xi\|_\infty<+\infty\},
\]
then $W$ is a Banach space with norm $\|\cdot\|_\infty$. By the
definition of $G_\lambda$, it is easy to check that there exists
$M>0$ that
\[
\|G_\lambda(\xi_1-\xi_2)\|_\infty\leq M\|\xi_1-\xi_2\|_\infty
\]
for any $\xi_1,\xi_2\in W$, i.e., ODE \eqref{equ 3.2 two ODE}
satisfies Lipschitz condition. As a result, according to the theory
of the linear ODE on the Banach space, ODE \eqref{equ 3.2 two ODE}
has the unique solution that
\[
F_t=e^{tG_\lambda}F_0
\]
for any $t\geq 0$. Since $F_0(x)=1$ for any $x\in \mathbb{Z}^d$,
\[
F_t(O)=\sum_{y:y\in Z^d}e^{tG_\lambda}(O,y)F_0(y)=\sum_{y:y\in
\mathbb{Z}^d}e^{tG_\lambda}(O,y).
\]
Since $G_\lambda(x,y)\geq 0$ when $x\neq y$,
$e^{tG_\lambda}(x,y)\geq 0$ for any $x,y\in \mathbb{Z}^d$.
Therefore, by Equation \eqref{equ 3.9},
\begin{equation}\label{equ 3.11}
E(\xi_t^2(O))=F_t(O)\leq \sum_{y\in
\mathbb{Z}^d}e^{tG_\lambda}(O,y)\frac{K(y)}{\inf_{x\in
\mathbb{Z}^d}K(x)}=\frac{K(O)}{\inf_{x\in \mathbb{Z}^d}K(x)}
\end{equation}
for any $t\geq 0$. According to Lemmas \ref{lemma 3.2}, \ref{lemma
3.3}, Equation \eqref{equ 3.11} and Cauchy-Schwartz inequality,
\begin{align}\label{equ 3.12}
\lim_{t\rightarrow+\infty}
P_\lambda\big(\eta_t(O)=1\big)&=\lim_{t\rightarrow+\infty}
P_\lambda\big(\xi_t(O)>0\big) \notag\\
&\geq
\limsup_{t\rightarrow+\infty}\frac{(E\xi_t(O))^2}{E(\xi_t^2(O))}
=\limsup_{t\rightarrow+\infty}\frac{1}{E(\xi_t^2(O))} \notag\\
&\geq \frac{\inf_{x\in \mathbb{Z}^d}K(x)}{K(O)}>0.
\end{align}
As a result,
\[
\lambda\geq \lambda_c
\]
for any $\lambda$ that there exists $K$ which satisfies $\inf_{x\in
\mathbb{Z}^d}K(x)>0$ and $G_\lambda K=0$.

\qed

\quad

\textbf{Acknowledgments.} The author is grateful to the financial
support from the National Natural Science Foundation of China with
grant number 11501542 and the financial support from Beijing
Jiaotong University with grant number KSRC16006536.

{}
\end{document}